\documentclass[10pt,reqno]{amsart}
\usepackage{amsmath}
\usepackage{amssymb}
\usepackage{amsfonts}
\makeatletter \@addtoreset{equation}{section} \makeatother

\newcommand{\noi}{\vspace{12pt}\noindent}
\newcommand{\beq}{\begin{equation}}
\newcommand{\eeq}{\end{equation}}
\newcommand{\bea}{\begin{eqnarray}}
\newcommand{\eea}{\end{eqnarray}}

\newcommand{\e}[1]{{(\ref{#1})}}
\newcommand{\eq}[1]{{eq.\ (\ref{#1})}}
\newcommand{\es}[2]{{(\ref{#1}) and (\ref{#2})}}
\newcommand{\eqs}[2]{{eqs.\ (\ref{#1}) and (\ref{#2})}}
\newcommand{\Ref}[1]{{Ref.~\cite{#1}}}
\newcommand{\mb}[1]{{\mbox{${#1}$}}}

\newcommand{\ie}{{${ i.e., \ }$}}
\newcommand{\eg}{{${ e.g., \ }$}}

\newcommand{\cf}{{cf.\ }}

\newcommand{\rhs}{{right--hand side }}

\renewcommand{\~}{ \ }

\newcommand{\eps}{\varepsilon^{}}
\newcommand{\cA}{\mathcal{A}}
\newcommand{\Z}{\mathbb{Z}}

\newcommand{\Ln}{{\rm Ln}}

\newcommand{\Hf}{\frac{1}{2}}

\newcommand{\papa}[1]{  \frac{\partial}{\partial {#1}}  }
\newtheorem{theorem}{Theorem}
\newtheorem{definition}[theorem]{Definition}

\newtheorem{proposition}[theorem]{Proposition}

\begin{document}
\thispagestyle{empty}
\title{Examples of Homotopy Lie Algebras}
\author{Klaus Bering}
\address{Institute for Theoretical Physics \& Astrophysics,
Masaryk University, Kotl\'a\v{r}sk\'a~2, CZ-611 37 Brno, Czech Republic}
\email{\tt bering@physics.muni.cz}
\author{Tom Lada}
\address{Department of Mathematics, North Carolina State University,
Raleigh NC 27695}
\email{\tt lada@math.ncsu.edu}

\begin{abstract}
We look at two examples of homotopy Lie algebras (also known as 
$L_{\infty}$ algebras) in detail from two points of view. We will 
exhibit the algebraic point of view in which the generalized Jacobi 
expressions are verified by using degree arguments and combinatorics. 
A second approach using the nilpotency of Grassmann-odd differential 
operators $\Delta$ to verify the homotopy Lie data is shown to produce 
the same results.
\end{abstract}

\maketitle

\section{Introduction}

\noi
Homotopy Lie algebras, or $L_\infty$ algebras, have been a topic 
of great interest to both mathematical physicists and to 
algebraists. By considering two different points of view, one can 
hope to gain a deeper understanding of these structures. In this 
note, we provide notations and definitions used by both communities, 
and hopefully illuminate both perspectives. On one hand, the second 
author and his collaborators \cite{daily,dailylada05} have
algebraicly constructed two concrete finite dimensional examples 
of homotopy Lie algebras from first principles. On the other hand, 
the first author and his collaborators have developed a
generalization of the Batalin-Vilkovisky formalism \cite{bv81} in 
which a nilpotent, Grassmann-odd, differential operator $\Delta$ 
may be used to identify $L_\infty$ structures, \cf Lemma in 
Sec.~2.3 of \Ref{bda96}, and Theorem~3.6 in \Ref{b06}. This
method is here applied to rederive the two examples of the 
second author and his collaborators.

\section{Homotopy Lie Algebras}
\label{sechla}

\noi
We begin by recalling the definition of an $L_\infty$ algebra 
\cite{ladastasheff93}, \cite{ladamarkl95}.

\begin{definition}
 An $L_\infty$ algebra structure on a $\Z$ graded vector space 
$V$ is a collection of graded skew symmetric linear maps
$l_n:V^{\otimes n}\rightarrow V$ of degree $2-n$ that satisfy 
generalized Jacobi identities
\beq
\sum_{i+j=n+1}\sum_{\sigma}e(\sigma)(-1)^{\sigma}
(-1)^{i(j-1)}l_j(l_i(
v_{\sigma(1)},\dots,v_{\sigma(i)}),
v_{\sigma(i+1)},\dots,v_{\sigma(n)})=0 \label{genjacid}
,
\eeq
where $(-1)^{\sigma}$ is the sign of the permutation, 
$e(\sigma)$ is the Koszul sign which is equal to -1 raised to the
product of  the degrees of 
the permuted elements, and $\sigma$ is taken 
over all $(i,n-i)$ unshuffles.
\label{def1}
\end{definition}

\noi
This is the cochain complex point of view; for chain complexes, 
require the maps $l_n$ to have degree $n-2$.

\subsection{Desuspension}

\noi
We will require an equivalent way to describe homotopy Lie 
algebra data that will be compatible with the operator approach.

\begin{definition}
Let $S^c(W)$ be the cofree cocommutative coassociative coalgebra 
on the graded vector space $W$. Then an $L_\infty$ algebra structure 
on $W$ is a coderivation $D:S^c(W)\rightarrow S^c(W)$ of degree $+1$ 
such that $D^2=0$.
\label{def2}
\end{definition}

\noi
Given an $L_\infty$ algebra structure $(V,l_i)$ as in 
Definition~\ref{def1}, we may desuspend $V$ to obtain the 
graded vector space $W=\downarrow V$, where $W_n=V_{n+1}$ and 
$\downarrow$ is the desuspension operator. Define 
$D:S^c(W)\rightarrow S^c(W)$ by 
$D=\hat{l}_0+\hat{l}_1+\hat{l}_2+\dots$,
where each $\hat{l}_n$ is a degree $+1$ symmetric map given by
\beq
\hat{l}_n = (-1)^{\frac{n(n-1)}{2}}\downarrow \circ l_n\circ 
\uparrow^{\otimes n}:S^c(W)\rightarrow W ,
\eeq
and then extended to a coderivation in the usual fashion. We will
demonstrate this construction explicitly in the examples.

\noi
The examples that we consider will be structures on relatively 
small graded vector spaces : $V=V_{0}\oplus V_{1}$, where each 
$V_{i}$ is finite dimensional. When we desuspend $V$, we will 
consider the graded vector space $W=W_{-1}\oplus W_{0}$.

\noi
We now describe the $\Delta$ operator approach.

\section{The $\Delta$ Operator Approach}
\label{secdelta}

\subsection{Vector Space $W$ with two Fermions}

\noi
To be concrete, we let $\dim(W_{-1})=2$. We use Greek indices 
$\alpha, \beta, \ldots \in \{1,2\}$ for a Fermionic basis
$\theta_{\alpha}\in W_{-1}$ with Grassmann parity 
$\eps(\theta_{\alpha})=1$. On the other hand, it will be useful to 
allow $W_{0}$ in the beginning to have infinitely many dimensions, 
and only at the very end perform a consistent truncation to a finite 
dimensional subspace. We use roman indices
$i,j,\ldots \in \{1,2,\ldots \}$ for the infinitely many Bosonic/even 
variables $x_{i}\in W_{0}$ with Grassmann parity $\eps(x_{i}) = 0$.
Hence, we are given a (super) vector space
\beq
W:=W_{-1}\oplus W_{0}, \quad 
W_{-1}:={\rm span}\langle\theta_{1},\theta_{2}\rangle , \quad 
W_{0}:={\rm span}\langle x_{1}, x_{2}, \ldots, x_{i}, \ldots \rangle .  
\eeq
We will for simplicity here only consider one kind of grading, 
although it is easy to generalize to several $\Z_2$ and $\Z$ 
gradings. In Section~\ref{sechla} we introduced a $\Z$ grading, 
called the degree. {}From an operational point of view, only a 
$\Z_2$ grading, the so-called Grassmann parity $\eps$, is needed. 
We shall start by only considering the $\Z_2$ grading $\eps$, and 
only later implement the full $\Z$ grading. This will lead to 
``selection rules'', \ie further restrictions.

\subsection{Algebra}

\noi
{}For an operational point of view, we use the fact the cocommutative 
coalgebra $S^c(W)$ has the same underlying vector space as the (super) 
symmetric algebra $\cA := {\rm Sym}^{\bullet}(W)$, where 
\beq
x_{i}\otimes x_{j} = x_{j}\otimes x_{i} , \qquad
x_{i}\otimes \theta_{\alpha}
=\theta_{\alpha}\otimes x_{i} , \qquad
\theta_{\alpha}\otimes\theta_{\beta}
=- \theta_{\beta}\otimes \theta_{\alpha} ,
\eeq
or
\beq 
z \otimes w = (-1)^{\eps(z)\eps(w)}w \otimes z
\eeq
for short, where $z,w \in W$.

\subsection{Bracket Hierarchy $\Phi^{\bullet}\equiv\hat{l}_{\bullet}$}

\noi
The family of maps $\hat{l}_{\bullet}$ on $W$ will be denoted by 
$\Phi^{\bullet}$ to conform with notation used in \Ref{bda96} and 
\Ref{b06}. We shall not always write (super) symmetric tensor 
symbol $\otimes$ explicitly. The sign convention is as follows:
\bea
\eps\left(\Phi^{n}(z_{1}\otimes z_{2}\otimes\ldots\otimes z_{n})\right) 
&=& 1+\eps(z_{1})+\eps(z_{2})+\ldots+\eps(z_{n}) , \cr
\Phi^{n}(\ldots\otimes z_{k} \otimes z_{k+1}\otimes\ldots)
&=&(-1)^{\eps(z_{k})\eps(z_{k+1})}
\Phi^{n}(\ldots\otimes z_{k+1}\otimes z_{k}\otimes\ldots) , \cr
\Phi^{n}(\lambda z_{1}\otimes z_{2}\otimes\ldots\otimes z_{n})
&=&(-1)^{\eps(\lambda)}\lambda
\Phi^{n}(z_{1}\otimes z_{2}\otimes\ldots\otimes z_{n}) ,\cr
\Phi^{n}(z_{1}\otimes\ldots\otimes z_{n}\lambda)
&=& \Phi^{n}(z_{1}\otimes\ldots\otimes z_{n})\lambda , \cr
z_{k}\lambda&=&(-1)^{\eps(z_{k})\eps(\lambda)} \lambda z_{k}~,
\eea 
Here $\lambda$ is a super number. We shall use multi-index notation
\beq
m = (m_{1}, m_{2}, \ldots, m_{i}, \ldots ) , \qquad
|m|=\sum_{i=1}^{\infty}m_{i} , \qquad
m!=\prod_{i=1}^{\infty}m_{i}! ,  \nonumber
\eeq
\beq
x^{\otimes m}=x_{1}^{\otimes m_{1}}\otimes x_{2}^{\otimes m_{2}}
\otimes\ldots\otimes x_{i}^{\otimes m_{i}}\otimes\ldots .
\eeq
The most general bracket hierarchy $\Phi^{\bullet}$ on $W$ is
\bea
\Phi^{|m|}(x^{\otimes m})
&=&c^{\alpha}_{m}\theta_{\alpha} , \label{phic} \\
\Phi^{|m|+1}(\theta_{\alpha}\otimes x^{\otimes m})
&=& b^{i}_{\alpha m}x_{i} , \label{phib} \\
\Phi^{|m|+2}(\theta_{\alpha}\otimes\theta_{\beta}\otimes x^{\otimes m})
&=& \epsilon_{\alpha\beta} a^{\gamma}_{m}\theta_{\gamma} , \label{phia}
\eea
where $a^{\alpha}_{m}$, $b^{i}_{\alpha m}$ and $c^{\gamma}_{m}$ are 
coefficients, and where 
\beq
\epsilon^{\alpha\beta} 
= -\epsilon^{\beta\alpha} , \qquad
\epsilon^{\alpha\beta}\epsilon_{\beta\gamma} 
= \delta^{\alpha}_{\gamma} , \qquad
\epsilon^{12} = 1 = \epsilon_{21} .
\eeq

\subsection{The $\Delta$ Operator}

\noi
Define generating functions
\beq
f^{\alpha}(p):=\sum_{m} a^{\alpha}_{m}\frac{p^{m}}{m!} , \qquad
g^{i}_{\alpha}(p):=\sum_{m}  b^{i}_{\alpha m}\frac{p^{m}}{m!} , \qquad
h^{\alpha}(p):=\sum_{m} c^{\alpha}_{m}\frac{p^{m}}{m!} .
\eeq
Define $\Delta$ operator 
\bea
\Delta &:=& \Delta_{2}+\Delta_{1}+\Delta_{0} ,\\
\Delta_{2}&:=&\Hf \theta_{\gamma} f^{\gamma}(\papa{x}) 
\epsilon_{\alpha\beta}
\papa{\theta_{\beta}}\papa{\theta_{\alpha}} , \\
\Delta_{1}&:=&x_{i} g^{i}_{\alpha}(\papa{x}) \papa{\theta_{\alpha}} , \\
\Delta_{0}&:=&\theta_{\alpha} h^{\alpha}(\papa{x}) .
\eea
We will from now on not always write the $\papa{x}$ dependence 
explicitly in the formula for $\Delta$.

\subsection{Koszul Bracket Hierarchy $\Phi_{\Delta}^{\bullet}$}

\noi
Define Koszul brackets hierarchy \mb{\Phi_{\Delta}^{\bullet}} as
\bea
\Phi_{\Delta}^{n}(z_{1}\otimes\ldots\otimes z_{n})
&:=&\underbrace{[[ \ldots [\Delta, L_{z_{1}}],\ldots ], L_{z_{n}}]}_{
n~{\rm commutators}}1 \label{phideltan} ,  \\
\Phi_{\Delta}^{0}&:=&\Delta(1)\equiv \theta_{\alpha} c^{\alpha}_{0} ,
\label{phidelta0}
\eea
where 
\beq
 L_{z}(w) := z w
\eeq
is the left multiplication operator with algebra element $z$. 

\noi
It is easy to check that the $\Phi_{\Delta}^{\bullet}$ Koszul brackets 
hierarchy \e{phideltan}-\e{phidelta0} reproduces the original 
$\Phi^{\bullet}$ bracket hierarchy \e{phic}-\e{phia}:
\beq
\Phi_{\Delta}^{\bullet}=\Phi^{\bullet} .
\eeq

\subsection{$L_{\infty}$ Structure and Nilpotency Conditions}

\noi
A consequence of Lemma in Sec.~2.3 of \Ref{bda96}, or alternatively 
Theorem~3.6 in \Ref{b06}, is that $\Phi_{\Delta}^{\bullet}$ forms a 
homotopy Lie algebra if and only if $\Delta$ is nilpotent (of 
order two), \ie $\Delta$ squares to zero,
\beq
\Delta^{2}\equiv \Hf [\Delta,\Delta]=0 .
\eeq
We calculate:
\bea
[\Delta_{2},\Delta_{2}]&=&0 ,\\
{}[\Delta_{2},\Delta_{1}]&=& \Hf x_{i}g^{i}_{\gamma}f^{\gamma}
\epsilon_{\alpha\beta}\papa{\theta_{\beta}}\papa{\theta_{\alpha}} ,
\label{d21}\\
{}[\Delta_{1},\Delta_{1}]
&=& 2x_{i} g^{i}_{\alpha}{}_{,j}g^{j}_{\beta}
\papa{\theta^{\alpha}}\papa{\theta^{\beta}}
= x_{i} g^{i}_{\alpha}{}_{,j} \epsilon^{\alpha\beta}g^{j}_{\beta}
\epsilon_{\gamma\delta}\papa{\theta^{\delta}}\papa{\theta^{\gamma}} , \\
{}[\Delta_{2},\Delta_{0}]
&=&\theta_{\gamma} f^{\gamma} h^{\alpha} \epsilon_{\alpha\beta}
\papa{\theta_{\beta}} , \\
{}[\Delta_{1},\Delta_{0}]&=&
\theta_{\alpha}h^{\alpha}{}_{,i}g^{i}_{\beta} \papa{\theta^{\beta}}
+x_{i} g^{i}_{\alpha} h^{\alpha} , \\
{}[\Delta_{0},\Delta_{0}]&=&0
\eea
{}For instance, \eq{d21} is proved as follows. Write shorthand 
$\Delta_{2}=\theta_{\gamma}D^{\gamma} $, where 
\beq
D^{\gamma} := \Hf f^{\gamma}(\papa{x})
\epsilon_{\alpha\beta}\papa{\theta_{\beta}}\papa{\theta_{\alpha}} .
\eeq
Then
\bea
[\Delta_{2},\Delta_{1}]&=&\theta_{\gamma}[D^{\gamma},\Delta_{1}]
+[\theta_{\gamma},\Delta_{1}]D^{\gamma}\cr 
&=&\theta_{\gamma}[D^{\gamma},\Delta_{1}]
+[\Delta_{1},\theta_{\gamma}]D^{\gamma} \cr
&=&\theta_{\gamma}[D^{\gamma},x_{i} g^{i}_{\delta} \papa{\theta_{\delta}}]
+[x_{i} g^{i}_{\alpha}\papa{\theta_{\alpha}},\theta_{\gamma}]D^{\gamma}\cr
&=&\theta_{\gamma}[D^{\gamma},x_{i}] g^{i}_{\delta} \papa{\theta_{\delta}}
+x_{i} g^{i}_{\alpha}[\papa{\theta_{\alpha}},\theta_{\gamma}]D^{\gamma}\cr
&=& \Hf \theta_{\gamma}f^{\gamma}{}_{,i} 
\epsilon_{\alpha\beta}\papa{\theta_{\beta}}\papa{\theta_{\alpha}}
g^{i}_{\delta} \papa{\theta_{\delta}}+x_{i} g^{i}_{\gamma}D^{\gamma} . 
\label{udregning}
\eea
Note that the first term on the \rhs of \e{udregning} must vanish 
because it contains three Fermionic derivatives, but there are only 
two different Fermions. The second term yields the result \e{d21}.

\noi
Altogether, the nilpotency condition $\Delta^{2}=0$ read
\bea
g^{i}_{\gamma}f^{\gamma}
+g^{i}_{\alpha}{}_{,j} \epsilon^{\alpha\beta}g^{j}_{\beta}&=&0 ,
\label{nilwronsk0} \\
f^{\alpha} h^{\gamma} \epsilon_{\gamma\beta}
+h^{\alpha}{}_{,i}g^{i}_{\beta} &=&0 , \\
g^{i}_{\alpha} h^{\alpha}&=&0 .\label{nilwronskh}
\eea

\subsection{Special Cases}

\noi
Let us now discuss special cases. Let us assume $h^{\alpha}\equiv 0$. 
Then the two last nilpotency conditions are satisfied, and only 
the first of the three nilpotency conditions 
\e{nilwronsk0}--\e{nilwronskh} 
remains.

\noi
Notice that we can explain $h^{\alpha}\equiv0$ as a ``selection rule'' 
from the degree $\Z$ grading, where $x_i\in W_{0}$ have degree $0$; 
$\theta_{\alpha}\in W_{-1}$ have degree $-1$; the brackets 
$\Phi_{\Delta}^{\bullet}$ have degree $+1$; and the $\Delta$ 
operator has degree $+1$. Then $c^{\alpha}_{m}\equiv0$,
$h^{\alpha}\equiv0$, and $\Delta_{0}\equiv0$.

\noi
Let us assume only one Bosonic/even variable $x\equiv x_{1}$,
\ie $0=x_{2}=x_{3}=\ldots$. Then the first nilpotency condition 
\e{nilwronsk0} reads:
\beq
g_{\gamma}f^{\gamma}+W(g_{1},g_{2})=0 ,  \label{nilwronsk}
\eeq
where
\beq
W(g_{1},g_{2}):=g^{\prime}_{\alpha} \epsilon^{\alpha\beta}g_{\beta}
\equiv g^{\prime}_{1} g^{}_{2}- g^{}_{1}g^{\prime}_{2}
\eeq
is the Wronskian.

\noi
Let us assume that $g_{1}$ is given with 
$g_{1}(p\!=\!0)\equiv b_{\alpha=1,m=0}\neq 0$.
Then we can interpret the inverse $1/g_{1}$ as a formal power series.

\noi
If there is also given $g_{2}$, then we can \eg choose
\beq
f^{1}=-\frac{W(g_{1},g_{2})}{g_{1}} , \label{nilwronsk2} \qquad
f^{2}=0 . 
\eeq
Or if there instead is also given $f^{1}$, then we can \eg choose
\beq
\frac{g_{2}}{g_{1}}= \int \! dp \~ \frac{f^{1}}{g_{1}} ,
\label{nilwronsk3} \qquad
f^{2}=0 . 
\eeq

\section{First Example}
\label{secex1}

\subsection{Algebra Approach}

\noi
The following $L_\infty$ algebra was studied in \cite{dailylada05}.
Let $V=V_0\oplus V_1$ be the graded vector space where $V_0$ has basis 
$\langle v_1,v_2\rangle$ and $V_1$ has basis $\langle w \rangle$. 
Define $l_n:V^{\otimes n}\rightarrow V$
 by
\beq
l_1(v_1)=l_1(v_2)=w , \qquad
l_2(v_1\otimes v_2)=v_1 , \qquad
l_2(v_1\otimes w)=w , \nonumber
\eeq
\beq
l_n\left(v_2\otimes w^{\otimes n-1}\right) =C_n w \mbox{ for } n \geq 3 ,
\eeq
and all other sectors are zero, and
where $C_n=(-1)^{\frac{(n-2)(n-3)}{2}}(n-3)!$. 

\noi
To verify the $L_\infty$ relations \e{genjacid}, the summands in 
the $L_{\infty}$ relation can be calculated as follows. The first summand 
reads
\beq
l_1 \circ l_n\left(v_1\otimes v_2\otimes w^{\otimes n-2}\right)=0 ,
\eeq
The next summand reads
\bea
\lefteqn{
l_2 \circ l_{n-1}\left(v_1\otimes v_2
\otimes w^{\otimes n-2}\right)} \cr
&=&(-1)^{n-1}l_2\left(l_{n-1}(
v_2\otimes w^{\otimes n-2})\otimes v_1\right) \cr
&=&(-1)^{n-1}C_{n-1}l_2(w\otimes v_1)=(-1)^n C_{n-1}w ,
\eea
{}For all \mb{3\leq k\leq n-3} we have 
\beq
l_k \circ l_{n-k+1}\left(v_1\otimes
v_2\otimes w^{\otimes n-2}\right)=0 ,
\eeq
because each summand in this expansion contains the term 
$l_k(v_1\otimes w^{\otimes k-1})=0$. The second-last summand reads
\bea
\lefteqn{
l_{n-1} \circ l_2\left(v_1\otimes v_2\otimes w^{\otimes n-2}\right)} \cr
&=&l_{n-1}\left(l_2(v_1\otimes v_2\otimes w^{\otimes n-2}\right) 
-(n-2)l_{n-1}\left(l_2(v_1\otimes w)
\otimes v_2\otimes w^{\otimes n-3}\right) \cr
&=&l_{n-1}(v_1\otimes w^{\otimes n-2}) 
-(n-2)l_{n-1}\left(w\otimes v_2\otimes w^{\otimes n-3}\right) \cr
&=&0+(n-2)l_{n-1}(v_2\otimes w^{\otimes n-2}) =(n-2)C_{n-2}w .
\eea
The last summand reads
\beq
l_n \circ l_1(v_1\otimes v_2\otimes w^{\otimes n-2})
=l_n(w\otimes v_2\otimes w^{\otimes n-2})
-l_n(w\otimes v_1\otimes w^{\otimes n-2})
=-C_n w .
\eeq
Consequently, the $n$th Jacobi expression is satisfied if and only if
\beq
\sum_{p=1}^n(-1)^{p(n-p)}
l_{n-p+1} \circ l_p\left(v_1\otimes v_2\otimes w^{\otimes n-2}\right)=0
\nonumber
\eeq
\beq
\Leftrightarrow 
(-1)^{(n-1)1}(-1)^n
C_{n-1}w+(-1)^{2(n-2)}(n-2)C_{n-1}w+(-1)^{1(n-1)}(-1)C_n w=0\nonumber
\eeq
\beq
\Leftrightarrow (-1)C_{n-1}+(n-2)C_{n-1}+(-1)^n C_n = 0
\Leftrightarrow C_n = (-1)^{n-1}(n-3)C_{n-1} .
\eeq
One can check that $C_n$ must equal $(-1)^{\frac{(n-2)(n-3)}{2}}(n-3)!$.

\subsection{Desuspension}

\noi
We next desuspend $V$ to obtain $W=W_{-1}\oplus W_0$ where $W_{-1}$ 
has basis $\langle \theta_1,\theta_2 \rangle$ and $W_0$ has basis 
$\langle x \rangle$ and then rewrite the $L_\infty$ data in terms 
of degree $+1$ maps
\beq 
\hat{l}_1(\theta_1)=\hat{l}_1(\theta_2)=x , \qquad
\hat{l}_2(\theta_1\otimes \theta_2)=\theta_1 , \nonumber
\eeq
\beq
\hat{l}_n\left(\theta_2\otimes x^{\otimes n-1}\right)=(-1)^n(n-3)!x ,
\eeq
and all other sectors are zero.
The $\hat{l}_n$'s will correspond to the $\Phi^n$'s in the next section.

\subsection{$\Delta$  Operator Approach}

\noi
The algebra of \Ref{dailylada05} has only one Bosonic generator 
$x\equiv x_{1}$, and is given as
\beq
\Phi^{2}(\theta_{1}\otimes\theta_{2})=\theta_{1} , \qquad
\Phi^{1}(\theta_{1})=x , \qquad
\Phi^{2}(\theta_{1}\otimes x)=x , \nonumber
\eeq
\beq
\Phi^{m+1}\left(\theta_2\otimes x^{\otimes m}\right)= 
\left\{\begin{array}{lcl} 
x &{\rm for}& m=0 \cr 
0 &{\rm for}& m=1 \cr
-(-1)^m (m-2)! \~ x &{\rm for}& m\geq 2 ,
\end{array} \right. 
\eeq
and all other sectors are zero. Thus the coefficients are
\beq
a^{1}_{m}= - \delta^0_m , \qquad
a^{2}_{m}= 0  , \qquad
b_{1m}= \delta^0_m+\delta^1_m , \nonumber
\eeq
\beq
b_{2m}= \left\{\begin{array}{lcl} 
1 &{\rm for}& m=0  \cr 
0 &{\rm for}& m=1  \cr
-(-1)^m (m-2)!&{\rm for}& m\geq 2 .
\end{array} \right.
\eeq
The generating functions become
\beq
f^1(p)= -1 ,\qquad
f^{2}(p)= 0 ,\qquad
g_1(p)= 1+p , \nonumber
\eeq
\beq
g_2(p)=1-\sum_{m=2}^{\infty}\frac{(-p)^m}{m(m-1)} 
= (1+p)[ 1- \ln(1+p)] .
\eeq
It is easy to check that the nilpotency condition \e{nilwronsk} is 
satisfied. Alternatively, $g_2$ could have been predicted from 
\eq{nilwronsk3}. 

\section{Second Example}
\label{secex2}

\subsection{Algebra Approach}

\noi
This next example was constructed by M.~Daily \cite{daily}.  Let 
$V=V_0\oplus V_1$ with $\dim(V_1) \geq \dim(V_0)$. Denote the basis 
for $V_0$ by $\langle v_1,...,v_i\rangle$ and the basis for $V_1$ 
by $\langle w_1,...,w_j \rangle$. Define
\beq
l_1(v_i)=w_i ,  \qquad
l_2(v_i\otimes v_j)=0 , \qquad
l_2(v_i\otimes w_j)=w_i+w_j , \nonumber
\eeq
\beq
l_n(v_i\otimes v_j\otimes w\mbox{-terms})=0 , \qquad
l_n(v_i\otimes w\mbox{-terms}) = C_n w_i \mbox{ for } n\geq 3 ,
\eeq
and all other sectors are zero. 

\noi
We begin verification of the $L_\infty$ algebra relations \e{genjacid} 
with
\bea
\lefteqn{
l_1\circ l_2(v_i\otimes v_j)-l_2\circ l_1(v_i\otimes v_j)} \cr
&=&l_1(0)-\left[l_2(l_1(v_i)\otimes v_j)-l_2(l_1(v_j)\otimes v_i)\right] 
\cr
&=&0-l_2(w_i\otimes v_j)+l_2(w_j\otimes v_i)=w_j+w_i-(w_i+w_j)=0 .
\eea
We next consider the generalized Jacobi expression evaluated 
on $v_i\otimes v_j\otimes w_k$. The first summand reads
\beq
l_1\circ l_3(v_i\otimes v_j\otimes w_k)=0 .
\eeq
The next summand reads
\bea
\lefteqn{
l_2\circ l_2(v_i\otimes v_j\otimes w_k)} \cr
&=&l_2(l_2(v_i\otimes v_j)\otimes w_k)
-l_2(l_2(v_i\otimes w_k)\otimes v_j) \cr
&&+l_2(l_2(v_j\otimes w_k)\otimes v_i) \cr
&=&0-l_2((w_i+w_k)\otimes v_j)+l_2((w_j+w_k)\otimes v_i) \cr
&=&(w_j+w_i)+(w_j+w_k)-(w_i+w_j)-(w_i+w_k) \cr 
&=&-w_i+w_j .
\eea
The last summand reads
\beq
l_3\circ l_1(v_i\otimes v_j\otimes w_k)
=l_3(w_i\otimes v_j\otimes w_k)
-l_3(w_j\otimes v_i\otimes w_k)
=-C_3w_j+C_3w_i.
\eeq
Thus, the generalized Jacobi expression
\beq
(l_1\circ l_3+l_2\circ l_2+l_3\circ l_1)
(v_i\otimes v_j\otimes w_k)=0 \nonumber
\eeq
\beq
\Leftrightarrow w_i+w_j-C_3w_j+C_3w_i=0
\Leftrightarrow C_3=1 .
\eeq
{}For $n\geq 4$, we compute
\beq
\sum_{p=1}^n(-1)^{p(n-p)}l_{n-p+1}\circ l_p(
v_i\otimes v_j\otimes w_{k_1}\otimes \dots \otimes w_{k_{n-2}}) .
\eeq
The first summand with $p=1$ reads
\bea
\lefteqn{
l_n\circ l_1(v_i\otimes v_j\otimes 
w_{k_1}\otimes \dots \otimes w_{k_{n-2}})} \cr
&=&l_n(w_i\otimes v_j\otimes w_{k_1}
\otimes \dots \otimes w_{k_{n-2}}) \cr
&&-l_n(w_j\otimes v_i\otimes w_{k_1}
\otimes \dots \otimes w_{k_{n-2}}) \cr
&=&-C_nw_j+C_nw_i 
=C_n(w_i-w_j) .
\eea
The next summand with $p=2$ reads
\bea
\lefteqn{
l_{n-1}\circ l_2(v_i\otimes v_j\otimes 
w_{k_1}\otimes \dots \otimes w_{k_{n-2}})} \cr 
&=&-\sum_{\alpha}l_{n-1}(l_2(v_i\otimes 
w_{k_{\alpha}})\otimes v_j\otimes w\mbox{-terms}) \cr
&&+\sum_{\alpha}l_{n-1}(l_2(v_j\otimes 
w_{k_{\alpha}})\otimes v_i\otimes w\mbox{-terms}) \cr
&=&-\sum_{\alpha}l_{n-1}((w_i+w_{k_{\alpha}})
\otimes v_j\otimes w\mbox{-terms}) \cr
&&+\sum_{\alpha}l_{n-1}((w_j+w_{k_{\alpha}})
\otimes v_i\otimes w\mbox{-terms}) \cr
&=&2(n-2)C_{n-1}(w_j-w_i) .
\eea
{}For $3\leq p\leq n-2$, we have
\bea
\lefteqn{
l_{n-p+1}\circ l_p(v_i\otimes v_j\otimes 
w_{k_1}\otimes \dots \otimes w_{k_{n-2}})} \cr
&=&(-1)^{p-1}\binom{n-2}{p-1}l_{n-p+1}\left(
l_p(v_i\otimes w-\mbox{terms})
\otimes v_j\otimes w\mbox{-terms}\right) \cr
&&-(-1)^{p-1}\binom{n-2}{p-1}l_{n-p+1}\left(l_p(v_j\otimes 
w\mbox{-terms})\otimes v_i\otimes w\mbox{-terms}\right) \cr
&=&(-1)^{p-1}\binom{n-2}{p-1}l_{n-p+1}\left(
C_pw_i\otimes v_j\otimes w\mbox{-terms}\right) \cr
&&-(-1)^{p-1}\binom{n-2}{p-1}l_{n-p+1}\left(
C_pw_j\otimes v_i\otimes w\mbox{-terms}\right) \cr
&=&(-1)^p\binom {n-2}{p-1}C_{n-p+1}C_pw_j
-(-1)^p\binom {n-2}{p-1}C_{n-p+1}C_pw_i \cr
&=&(-1)^{p+1}\binom {n-2}{p-1}C_{n-p+1}C_p(w_i-w_j) .
\eea
The second-last summand with $p=n-1$ reads
\bea
\lefteqn{
l_2\circ l_{n-1}(v_i\otimes v_j\otimes w_{k_1}
\otimes \dots \otimes w_{k_{n-2}})} \cr
&=&(-1)^{n-2}l_2(l_{n-1}(v_i\otimes w\mbox{-terms})\otimes v_j) 
-(-1)^{n-2}l_2(l_{n-1}(v_j\otimes w\mbox{-terms})\otimes v_i) \cr
&=&(-1)^{n-2}l_2(C_{n-1}w_i\otimes v_j)
-(-1)^{n-2}l_2(C_{n-1}w_j\otimes v_i) \cr
&=&(-1)^{n-1}C_{n-1}(w_j+w_i)
-(-1)^{n-1}C_{n-1}(w_i+w_j)
=0 .
\eea
The last summand with $p=n$ reads
\beq
 l_1\circ l_n(v_i\otimes v_j\otimes 
w_{k_1}\otimes \dots \otimes w_{k_{n-2}})=0 .
\eeq
We add together all of the above summands with $p=1,2,\ldots,n$ to obtain
\bea
\lefteqn{
\sum_{p=1}^n(-1)^{p(n-p)}l_{n-p+1}\circ 
l_p(v_i\otimes v_j\otimes w_{k_1}
\otimes \dots \otimes w_{k_{n-2}})} \cr
&=&(-1)^{n-1}C_n(w_i-w_j)-2(n-2)C_{n-1}(w_i-w_j) \cr
&&+\sum_{p=3}^{n-2}(-1)^{p(n-p)}(-1)^{p+1}
\binom{n-2}{p-1}C_{n-p+1}C_p(w_i-w_j)+0+0.
\eea
So,
\beq
\sum_{p=1}^n(-1)^{p(n-p)}l_{n-p+1}\circ 
l_p(v_i\otimes v_j\otimes w_{k_1}
\otimes \dots \otimes w_{k_{n-2}})=0 \nonumber
\eeq
\beq
\Leftrightarrow (-1)^{n-1}C_n-2(n-2)C_{n-1}
+\sum_{p=3}^{n-2}(-1)^{pn+1}\binom{n-2}{p-1}C_{n-p+1}C_p=0.
\eeq
One can then solve for
\beq
C_n = (-1)^n\left[-2(n-2)C_{n-1}
+\sum_{p=3}^{n-2}(-1)^{pn+1} \binom{n-2}{p-1} C_{n-p+1}C_p\right]
\eeq
with $C_3=1$.  

\subsection{Desuspension}

\noi
As before, we desuspend the vector space to obtain 
$W=W_{-1}\oplus W_0$ and convert the $l_n$'s to degree $+1$ symmetric 
maps and end up with the homotopy Lie algebra structure given by
\beq
\hat{l}_1(\theta_i) = x_i , \qquad
\hat{l}_2(\theta_i\otimes \theta_j)=0 , \qquad
\hat{l}_2(\theta_i\otimes x_j)=x_i+x_j , \nonumber
\eeq
\beq
\hat{l}_n(\theta_i\otimes \theta_j\otimes x\mbox{-terms})=0 , \nonumber
\eeq
\beq
\hat{l}_n(\theta_i\otimes x\mbox{-terms})
=(-1)^{\frac{n(n-1)}{2}}(-1)^{n-1}C_n x_i \mbox{ for }  n\geq 3 ,
\eeq
and all other sectors are zero. The last equation may be rewritten as
\beq
\hat{l}_n(\theta_i\otimes x\mbox{-terms})
=(2-n)^{n-2}x_i \mbox{ for }  n\geq 3 . 
\eeq

\subsection{$\Delta$ Operator Approach}

\noi
In the following we let $\dim(W_{-1})=2$, to conform with the theory
developed in Section~\ref{secdelta}. Moreover, it is practical to
let $W_{0}$ have infinitely many Bosonic generators $x_i$. (It will be 
consistent to truncate the tail $0=x_{N+1}=x_{N+2}=\ldots$ to reduce to 
only finitely many generators $x_1, \ldots, x_{N}$.) Then the second 
example is of the form
\bea
\Phi^1(\theta_{\alpha})&=&B_0 x_{\alpha} , \cr
\Phi^2(\theta_{\alpha}\otimes x_{i})&=&B_1 x_{\alpha}+x_{i} ,\cr
\Phi^{|m|+1}(\theta_{\alpha}\otimes x^{\otimes m})
&=&B_{|m|}x_{\alpha} \mbox{ for } |m| \geq 2 , \label{secondex}
\eea
and all other sectors are zero, and where $B_0, B_1, B_2, \ldots$ are 
complex numbers with \mb{B_0\neq 0}. By scaling 
\beq
x'_{i}=B_0 x_{i} , \quad 
\theta'_{\alpha}=\theta_{\alpha} , \quad 
\Phi^{'}=\Phi , \quad 
B'_{M}=(B_0)^{M-1}B_{M}\mbox{ for } M=0,1,2, \ldots ,
\eeq
(and by dropping the primes again afterwards) we will from now on always 
assume the initial condition 
\beq
B_{0}=1 .  \label{ic1}
\eeq
We will below prove the following Proposition~\ref{prop3}.

\begin{proposition} 
The $\Phi^{\bullet}$ bracket hierarchy \e{secondex} with 
initial condition \e{ic1} is a homotopy Lie algebra if and only if 
\beq
B_{M} = (1-M)^{M-1}\mbox{ for } M=0,1,2,\ldots \label{bcoeff}
\eeq
(with the convention that $0^0:=1$).
\label{prop3}
\end{proposition}

\noi
{\bf Proof}. The bracket coefficients are in this example
\bea
a^{\alpha}_{m}&=& 0, \cr
b^i_{\alpha m}
&=&\delta^i_{\alpha}B_{|m|}+\delta^0_{m_1}\delta^0_{m_2}\cdots
\delta^0_{m_{i-1}}\delta^1_{m_i}\delta^0_{m_{i+1}}\cdots . 
\eea
The generating functions become
\bea
f^{\alpha}(p)&=& 0 , \cr
g^i_{\alpha}(p)&=&\delta^i_{\alpha}G(P)+p^i ,
\eea
where
\beq
G(P)=\sum_{m} B_{|m|}\frac{p^m}{m!} 
= \sum_{M=0}^{\infty}B_{M}\frac{P^M}{M!} ,
\label{gexpan}
\eeq
and
\beq
P:=\sum_{i=1}^{\infty}p^i .
\eeq
The initial condition \e{ic1} becomes
\beq
G(P\!=\!0)=1 . \label{ic2}
\eeq
The nilpotency condition \e{nilwronsk0} reads
\bea
(\alpha \leftrightarrow \beta) &=&g^{i}_{\alpha}{}_{,j}g^{j}_{\beta} 
=(\delta^i_{\alpha}G'(P)+\delta^i_j)(\delta^j_{\beta}G(P)+p^j) \cr
&=&\delta^i_{\alpha}G'(P)(G(P)+P)+\delta^i_{\beta}G(P)+p^i .
\eea
This is equivalent to the ODE
\beq
G'(P)(G(P)+P)=G(P) \label{gode}
\eeq
\beq 
\Leftrightarrow \frac{dP}{dG}=1+\frac{P}{G} 
\Leftrightarrow \frac{d}{dG}\left[\frac{P}{G}\right]=\frac{1}{G} 
\Leftrightarrow \frac{P}{G}=\Ln(G)+{\rm constant} .
\eeq
We deduce from the initial condition \e{ic2} that the inverse function 
$P=P(G)$ is
\beq
P(G)=G\Ln(G) = -(1-G)+\sum_{n=2}^{\infty}\frac{(1-G)^n}{n(n-1)} .
\label{invg}
\eeq
Let us now recall the Lambert function $W=W(P)$, whose inverse function 
$P=P(W)$ is
\beq
P(W)=We^{W} \label{invw} .
\eeq
(Hopefully, the reader will not be confused by the fact that we denote 
two different function $P=P(G)$ and $P=P(W)$ (and in fact also the 
``momentum'' variable $P$ itself) with the same symbol $P$. It should be 
clear from the context which is which.) Note that the Lambert function 
$W=W(P)$ has a zero in $P=0$ 
\beq
W(P\!=\!0)=0 .  \label{ic3}
\eeq 
By comparing \eqs{invg}{invw} we deduce that the sought-for function 
$G=G(P)$ is just the exponential of the Lambert function
\beq
G(P)=e^{W(P)}=\frac{P}{W(P)} .
\eeq
The Taylor expansion for the Lambert function $W=W(P)$ is 
\beq
W(P)= \sum_{n=1}^{\infty} (-n)^{n-1} \frac{P^n}{n!} . \label{taylorw}
\eeq
The Taylor coefficients with $n\geq 1$ follow from Lagrange's inversion 
formula, or simply by calculating
\bea
W^{(n)}(P\!=\!0)&=& \frac{1}{n!} \oint_{0} \frac{dP}{2\pi 
i}\frac{W'(P)}{P^n}
= \frac{1}{n!} \oint_{0} \frac{dW}{2\pi i}\frac{1}{P(W)^n} \cr
&=& \frac{1}{n!} \oint_{0} \frac{dW}{2\pi i}\frac{e^{-nW}}{W^n}
=\left. \frac{d^{n-1}}{dW^{n-1}} e^{-nW} \right|_{W\!=\!0}\cr
&=&(-n)^{n-1} .
\eea
Similarly, the Taylor coefficients for the function $G=G(P)$ with 
$n\geq 1$ are
\bea
B_{n}&=&G^{(n)}(P\!=\!0)
= \frac{1}{n!} \oint_{0} \frac{dP}{2\pi i}\frac{G'(P)}{P^n}
= \frac{1}{n!} \oint_{0} \frac{dP}{2\pi i}\frac{W'(P) e^{W(P)}}{P^n} \cr
&=& \frac{1}{n!} \oint_{0} \frac{dW}{2\pi i}\frac{e^W}{P(W)^n}
= \frac{1}{n!} \oint_{0} \frac{dW}{2\pi i}\frac{e^{(1-n)W}}{W^n} 
= \left. \frac{d^{n-1}}{dW^{n-1}} e^{(1-n)W} \right|_{W\!=\!0} \cr
&=&(1-n)^{n-1} .
\eea
The Taylor expansion for the function $G=G(P)$ is 
\beq
G(P)= \sum_{n=0}^{\infty} (1-n)^{n-1} \frac{P^n}{n!} . \label{taylorg}
\eeq
Both Taylor series \es{taylorw}{taylorg} have radius of convergence 
equal to $1/e$, as may be seen by the ratio test. This completes the 
proof of Proposition~\ref{prop3}.

\vspace{0.8cm}

\noi
{\sc Acknowledgement:}~The work of K.B.\ is supported by the Ministry
of Education of the Czech Republic under the project MSM 0021622409.


\begin{thebibliography}{999}

\bibitem{bv81} 
I.A.~Batalin and G.A.~Vilkovisky, {\em Gauge algebra and quantization}, 
Phys.~Lett.\ {\bf 102B} (1981) 27-31.

\bibitem{bda96}
K.~Bering, P.H.~Damgaard and J.~Alfaro, {\em Algebra of higher antibrackets}, 
Nucl.~Phys.\ {\bf B478} (1996) 
459-504.

\bibitem{b06} 
K.~Bering, {\em Non-commutative Batalin-Vilkovisky algebras, 
homotopy Lie algebras and the Courant bracket}, 
Commun.~Math.~Phys.\ {\bf 274} (2007) 297-341. 

\bibitem{daily}
M.~Daily, {\em Examples of $L_m$ and $L_\infty$ structures on $V_0\oplus V_1$},
unpublished notes.

\bibitem{dailylada05} M.~Daily and T.~Lada, 
{\em A finite dimensional $L_\infty$ algebra example in gauge theory}, 
Homotopy, Homology and Applications Vol.\ {\bf 7} (2005) 87-93.

\bibitem{ladamarkl95}
T.~Lada and M.~Markl, {\em Strongly homotopy Lie algebras}, 
Commun.~Algebra {\bf 23} (1995) 2147-2161.

\bibitem{ladastasheff93} T.~Lada and J.D.~Stasheff, 
{\em Introduction to sh Lie algebras for physicists}, 
Int.~J.~Theor.~Phys.\ {\bf 32} (1993) 1087-1103.


\end{thebibliography}
\end{document}